\documentclass[11pt]{amsart}
\usepackage{geometry}                
\usepackage{graphicx}
\usepackage{ctable}
\usepackage{booktabs}
\usepackage{tabularx}
\usepackage{multirow}
\usepackage{bigints}
\usepackage{amssymb}
\usepackage{epstopdf}
\usepackage[justification=raggedright, singlelinecheck=false]{caption}
\usepackage{subfig}
\usepackage{rccol}
\usepackage{verbatim}
\usepackage{amsthm}
\usepackage{float}
\usepackage{nccmath}
\usepackage{pinlabel}
\newcommand{\p}{\\[10 pt]}
\newtheorem{conjec}{Conjecture}

\newtheorem{prop}{Proposition}

\title{iterations of quadratic polynomials over finite fields}
\author{William Worden}

\begin{document}

\begin{abstract}
Given a map \(\mathit{f}:\mathbb{Z}\rightarrow\mathbb{Z}\) and an initial argument \(\alpha\), we can iterate the map to get a finite set of iterates modulo a prime \(p\).  In particular, for a quadratic map \(\mathit{f}(z)=z^2 +c\), \(c\) constant, work by Pollard suggests that this set should have length on the order of \(\sqrt{p}\).  We give a heuristic argument that suggests that the statistical properties of this set might be very similar to the Birthday Problem random variable \(X_n\), for an \(n=p\) day year, and offer considerable experimental evidence that the limiting distribution of these set lengths, divided by \(\sqrt{p}\), for \(p\leq x\) as \(x \rightarrow \infty\), converges to the limiting distribution of \(X_n/\sqrt{n}\), as \(n\rightarrow \infty\).
\end{abstract}

\maketitle

\section{introduction}
Let \(\mathit{f}\in\mathbb{Z}[z]\) be a polynomial and let \(\alpha\in\mathbb{Z}\).  We define the orbit of \(\alpha\) under \(\mathit{f}\) to be
\[
\mathcal{O}_\mathit{f}(\alpha)=\{\mathit{f}^n(\alpha) : n=0,1,2,3,\dots\},
\]
and for each prime \(p\) we define the orbit modulo \(p\) of \(\alpha\) under \(\mathit{f}\) to be
\[
\mathcal{O}_\mathit{f}^p(\alpha)=\{\mathit{f}^n(\alpha)\bmod p : n=0,1,2,3,\dots\},
\]
where \(\mathit{f}^n\) is the \(n^{th}\) iterate of \(\mathit{f}\),
\[
\mathit{f}^n=\underbrace{\mathit{f}\circ\mathit{f}\circ\dots\circ\mathit{f}}_n,
\]
and \(\mathit{f}^0(\alpha)=\alpha\). For a fixed \(\mathit{f}\) and \(\alpha\) and a given prime \(p\), let \(m_p\) be the order of \(\mathcal{O}_\mathit{f}^p(\alpha)\).

If $f$ is a random map, i.e., a map chosen from the uniformly distributed set consisting of all maps from $\mathbb{F}_p$ into $\mathbb{F}_p$ (see Harris \cite{harris}), then the values of \(\mathit{f}^n(\alpha) \) are uniformly distributed for all $n$, and all $\alpha$, and so the probability that \(\mathit{f}^0(\alpha),\mathit{f}^1(\alpha),\mathit{f}^2(\alpha),...,\mathit{f}^k(\alpha)\) are all different is
\[
1\cdot\frac{p-1}{p}\cdot\frac{p-2}{p}\cdot\dots\cdot\frac{p-k}{p}=\frac{(p-1)!}{p^k(p-k-1)!},
\]
since, once \(\alpha\) is fixed, there are \(p-1\) choices for \(\mathit{f}^1(\alpha)\),  \(p-2\) choices for \(\mathit{f}^2(\alpha)\), and so on.  Therefore, in this case the probability that (at least) two of \(\mathit{f}^0(\alpha),\mathit{f}^1(\alpha),\mathit{f}^2(\alpha),\dots,\mathit{f}^k(\alpha)\) are equal is
\[
q_k^{(p)}=1-\frac{(p-1)!}{p^k(p-k-1)!}.
\]
By an analogous argument, \(q_k^{(p)}\) is also the probability that among \(k\) people, two people have the same birthday, where \(p\) is the number of days in a year.  Framing this a little differently, we let the random variable \(X_n\) be the number of times that we must sample (uniformly, with replacement) from the set \(\{1,2,3,\dots,n\}\) to get a repetition.  Since it is known that the expected value of this variable is on the order of \(\sqrt{n}\), we look instead at the variable \(\frac{X_n}{\sqrt{n}}\). 

In light of the above heuristic, we might expect that for a fixed polynomial $\mathit{f}$ and initial value $\alpha$, $\frac{m_p}{\sqrt{p}}$ will, on average, ``behave" similarly to  \(\frac{X_n}{\sqrt{n}}\).  In particular, we might guess that the limiting distribution of $\frac{m_p}{\sqrt{p}}$, for $p \leq x$, \(x\rightarrow\infty$, will be similar to the limiting distribution of \(\frac{X_n}{\sqrt{n}}\), as $n \rightarrow \infty$.  We note that the above heuristic is not new; similar arguments have been given by Pollard \cite{pollard}, Bach \cite{bach}, and Brent \cite{brent} to name a few, leading to conjectures that $m_p$ is on average $\approx \sqrt{\frac{\pi}{2}p}$.

We also consider a related question.  For a fixed \(\mathit{f}\in\mathbb{Z}[z]\), \(\alpha\in\mathbb{Z}\), let 
\[
\mathcal{Q}_{\mathit{f},\alpha}(x)=\{p\leq x : \mathit{f}^n(\alpha)\equiv 0 \!\!\!\pmod p \quad\text{for some } n=0,1,2,\dots\}.
\]
That is, \(\mathcal{Q}_{\mathit{f},\alpha}(x)\) is the set of primes \(p\) less than or equal to \(x\) such that \(0\) appears in the orbit modulo \(p\) of \(\alpha\) under \(\mathit{f}\).  In particular, we are interested in the size of  \(\mathcal{Q}_{\mathit{f},\alpha}(x)\). Since, for a given prime $p$, the proportion of elements mod $p$ in the orbit of $\alpha$ under $\mathit{f}$ is $\frac{m_p}{p}$, we hypothesize that \(|\mathcal{Q}_{\mathit{f},\alpha}(x)|\) will grow at a rate proportional to \(\frac{m_p}{p}\). Therefore, if we are correct that \(m_p\) will grow at a rate proportional to \(\sqrt{p}\), we might expect that
\[
|\mathcal{Q}_{\mathit{f},\alpha}(x)|= \sum_{p\leq x} \frac{m_p}{p}=c \cdot  \frac{\sqrt{x}}{\log x},
\]
for some constant \(c\in\mathbb{R}\).  

In the following we take an experimental approach to studying properties of the set $\frac{m_p}{\sqrt{p}}$.  For selected maps \(\mathit{f}\) and initial values \(\alpha\), we compute the orbits modulo \(p\) for all \(p\leq 2^{25}\).  In particular, given these orbits we can find the moments of $\frac{m_p}{\sqrt{p}}$, and the length of \(\mathcal{Q}_{\mathit{f},\alpha}(x)\). As we will demonstrate in the sections to follow, our results give strong support to the above heuristic,  and lead us to make the following conjectures:

\begin{fleqn}
\begin{conjec}
Let \(\mathit{f}(z)=z^2 + c\) and \(\alpha\in\mathbb{Z}\) be such that
\[\quad 1) \quad c\in \mathbb{Z}\setminus\{0,-2\}\]

\[\quad 2) \quad\alpha \neq \scriptstyle\pm\frac{1}{2}\displaystyle(1\pm\scriptstyle\sqrt{1-4c}\displaystyle), \quad \alpha \neq \scriptstyle\pm \frac{1}{2}\displaystyle (1 \pm \scriptstyle\sqrt{-3-4c}\displaystyle), \quad \alpha \neq 0,\pm1 \text{ when } c=-1,\]
and let the orbit length \(m_p\) be as defined above.  Then as \( x \rightarrow \infty$, the distribution of $\frac{m_p}{\sqrt{p}}$ converges, independent of $\mathit{f}$ and $\alpha$, to a continuous distribution \(F(t)=1-e^{-t^2/2}, t\geq 0.\)
In particular, the \(r^{th}\) moments  of $\frac{m_p}{\sqrt{p}}$ are given by \(\mu_r=r(r-2)(r-4)\cdots 2\) for \(r\) even, and \(\mu_r=r(r-2)(r-4)\cdots 1\cdot \sqrt{\frac{\pi}{2}}\) for \(r\) odd.
\end{conjec}
\end{fleqn}
The motivation for the result conjectured above is elaborated upon in section 2, and the need to include conditions (1) and (2) for both conjectures is explained in section 4.
\p
\begin{conjec}
Let \(\mathit{f}(z)=z^2+c\) and \(\alpha\in\mathbb{Z}\) be such that conditions \((1)\) and \((2)\) of Conjecture 1 hold, and \(\alpha^2 \neq -c\). Define \(\mathcal{Q}_{\mathit{f},\alpha}(x)=\{p\leq x : \mathit{f}^n(\alpha)\equiv 0 \pmod p \,\,\,\text{for some } n\geq0\}\). Then
\[
\lim_{x\rightarrow\infty} \left(|\mathcal{Q}_{\mathit{f},\alpha}(x)| \frac{\log x}{\sqrt{x}}\right)=\sqrt{2\pi}
\]
\end{conjec}

\subsection*{Acknowledgements}
Our research was funded by a Rich Summer Internship grant from the Dr. Barnett and Jean-Hollander Rich Scholarship Fund. We are very grateful for the opportunity that this grant allowed, and would like to thank the selection committee and those who have made contributions to the fund. We would also like to thank our advisor Gautam Chinta, whose guidance and assistance throughout the research process were instrumental to the project's success, and whose invaluable feedback during the drafting of this paper improved its quality considerably. In addition, we thank Prof. Hutz for calling our attention to the paper \emph{Periods of Rational Maps Modulo Primes}, wherein the authors carry out similar computations and obtain results compatible with our computations, and we thank Prof. Silverman for directing us to the publication in which his article, \emph{Variation of Periods Modulo $p$ in Arithmetic Dynamics}, appeared. \cite{hutz, silv}.

\section{length of the orbit modulo p and the birthday problem}
Let \(E_k\) be the \(k^{th}\) number drawn uniformly from the set \(\{1,2,3,\dots,n\}\), with replacement, and let \(X_n\) be as defined in section 1. Then for \(k\leq n\) we have

\begin{align*}
P(X_n>k)&=P(E_1,\dots,E_k \,\,\text{all take different values})\\
&=\prod_{j=2}^{k} \big(1-P(E_j=E_i \,\,\text{for some} \,\,i<j)\big)\\
&=\prod_{j=2}^{k} \bigg( 1-\frac{j-1}{n}\bigg)=\exp{\bigg[\sum_{j=1}^{k-1} \log{(1-j/n)}\bigg]}
\end{align*}

So as \(n\rightarrow\infty\), we have the following for \(0 \leq t \leq \sqrt{n}\):

\begin{align*}
\lim_{n\rightarrow \infty} P(X_n/\sqrt{n} > t)&=\lim_{n\rightarrow \infty} P(X_n> t\sqrt{n})\\
&=\lim_{n\rightarrow \infty} \exp{\bigg[\sum_{1\leq j < t\sqrt{n}} \log{(1-j/n)}\bigg]}\\
&=\lim_{n\rightarrow \infty} \exp{\bigg[-\sum_{1\leq j < t\sqrt{n}} \sum_{k=1}^{\infty} \frac{(j/n)^k}{k} \bigg]}
\end{align*}
where we have used the power series representation for $\log{(1-j/n)}$ in the third line. Switching the order of summation, and pulling the first term of the sum over $k$ out of the exponential, we have
\begin{align*}
&=\lim_{n\rightarrow \infty} \exp{\bigg[-\sum_{1\leq j < t\sqrt{n}} j/n\bigg]}\cdot \lim_{n\rightarrow \infty} \exp{\bigg[-\sum_{k=2}^{\infty} \sum_{1\leq j < t\sqrt{n}} \frac{(j/n)^k}{k}\bigg]}\\
&\approx \lim_{n\rightarrow \infty} \exp{\bigg[-\frac{t\sqrt{n}(t\sqrt{n}+1)}{2n}\bigg]}\cdot \lim_{n\rightarrow \infty}  \exp{\bigg[-\sum_{k=1}^{\infty} O\bigg(\frac{t^{k+2}}{kn^{k/2}}\bigg)\bigg]}\\
&\approx e^{-t^2/2}\cdot \exp{\bigg[\sum_{k=1}^{\infty} \lim_{n\rightarrow \infty} O\bigg(\frac{t^{k+2}}{kn^{k/2}}\bigg)\bigg]} = e^{-t^2/2}.
\end{align*}
where the second line follows because, in general, $\sum_{1\leq j\leq m} j^k$ is a polynomial in $m$ of degree $k+1$, and the third line, where we have brought the limit inside the sum, follows from the Monotone Convergence Theorem. Therefore \(\lim_{n\rightarrow\infty} P(X_n/\sqrt{n} \leq t)=1-e^{-t^2/2}\), so we see that the distribution of \(\frac{X_n}{\sqrt{n}}\) converges to a distribution function \(F(t)=1-e^{-t^2/2}\), which has an associated density function \(\mathit{f}(t)=F'(t)=te^{-t^2/2}\). To support our conjecture in section 1---that \(F(t)\) is the limiting distribution of $\frac{m_p}{\sqrt{p}}$, as \(x\rightarrow \infty\)---we compare the moments of $\frac{m_p}{\sqrt{p}}$, which we compute in section 5 for large \(x\), to the limiting moments of \(\frac{X_n}{\sqrt{n}}\), as \(n\rightarrow \infty\).  With the limiting density function \(\mathit{f}(t)\) of \(\frac{X_n}{\sqrt{n}}\) in hand we can derive a general expression for the \(r^{th}\) moment:

\begin{align*}
\mu_r=\int_0^{\infty} t^r\mathit{f}(t)\, \mathrm{d}t&=\int_0^{\infty} t^{r+1}e^{-t^2/2}\,\,\mathrm{d}t\\
&=\left.-t^re^{-t^2/2}\right |_0^{\infty} + r \int_0^{\infty} t^{r-1}e^{-t^2/2}\,\, \mathrm{d}t\\
&=r \int_0^{\infty} t^{r-1}e^{-t^2/2}\,\, \mathrm{d}t,
\end{align*}
where \(r\) applications of L'Hospital's rule give us \(0\) for the \(-t^re^{-t^2/2}\) term.  We continue applying integration by parts as above until we get
\begin{align*}
\mu_r &= r(r-2)(r-4)\cdots 2\cdot \int_0^{\infty} te^{-t^2/2}\,\,\mathrm{d}t\qquad\text{if \(r\) is even,}\\
\mu_r &= r(r-2)(r-4)\cdots 1\cdot \int_0^{\infty} e^{-t^2/2}\,\,\mathrm{d}t  \,\, \qquad\text{if \(r\) is odd,}
\end{align*}
The first integral above evaluates to \(\left.-e^{-t^2/2}\right|_0^{\infty}=1\), and the second integral we evaluate as follows:
\begin{align*}
I &= \int_0^{\infty} e^{-t^2/2}\,\,\mathrm{d}t\\
\implies \quad (2I)^2 &=\left(\int_{-\infty}^{\infty} e^{-t^2/2}\,\,\mathrm{d}t\right)^2\\
&=\int_{-\infty}^{\infty} e^{-x^2/2}\,\,\mathrm{d}x \cdot \int_{-\infty}^{\infty} e^{-y^2/2}\,\,\mathrm{d}y\\
&=\int_{-\infty}^{\infty} \int_{-\infty}^{\infty} e^{-(x^2+y^2)/2}\,\,\mathrm{d}x\,\mathrm{d}y\\
&=\int_{r=0}^{\infty}\int_{\theta=0}^{2\pi} re^{-r^2/2}\,\,\mathrm{d}r\,\mathrm{d}\theta = 2\pi\\
\implies \,\,\,\,\quad\quad I &=\sqrt{\frac{\pi}{2}}
\end{align*}
Therefore the \(r^{th}\) moments of the limiting distribution of \(\frac{X_n}{\sqrt{n}}\), as \(n \rightarrow \infty\), are given by
\begin{align*}
\mu_r &= r(r-2)(r-4)\cdots 2\qquad\,\,\,\,\quad\quad\quad\text{if \(r\) is even,}\\
\mu_r &= r(r-2)(r-4)\cdots 1\cdot \sqrt{\pi/2}  \,\, \qquad\text{if \(r\) is odd.}
\end{align*}
For the first four moments this gives us \(\mu_1=\sqrt{\pi/2}, \,\,\mu_2=2, \,\,\mu_3=3\sqrt{\pi/2}, \,\,\mu_4=8\).  Therefore, to support our claim in Conjecture 1 we must provide evidence that the moments of $\frac{m_p}{\sqrt{p}}$ are converging, as \(x\rightarrow \infty\), to the moments \(\mu_r\) above.  In our computations we use the following expression for the \(r^{th}\) moments of $\frac{m_p}{\sqrt{p}}$:
\[
M_r=\frac{1}{|\{p\leq x\}|}\sum_{p\leq x} \left(\frac{m_p}{\sqrt{p}}\right)^r.
\]

\section{iterates of \(\mathit{f}\) congruent to zero modulo p}
In this section we consider the quantity \(|\mathcal{Q}_{\mathit{f},\alpha}(x)|\frac{\log x}{\sqrt{x}}\), as defined in Section 1.  Assuming that the probability that \( 0 \in \mathcal{O}_\mathit{f}^p(\alpha)\) is \(\frac{m_p}{p}\),  and that \(M_1\) will converge to \(\sqrt{\scriptstyle \pi/2}\), we define \(G(x)=\frac{\log x}{\sqrt{x}} \sum_{p\leq x} \frac{\sqrt{\scriptstyle \pi/2}}{\sqrt{p}}\), and make a guess that
\begin{equation}
\label{guess0}
\lim_{x\rightarrow\infty} \left(|\mathcal{Q}_{\mathit{f},\alpha}(x)| \frac{\log x}{\sqrt{x}}\right)= \lim_{x\rightarrow\infty} G(x).
\end{equation}
If we let \(\pi(x)=\sum_{k\leq x} a(k)\), where \(a(k)=1\) if \(k\) is prime and \(0\) otherwise, and define \(\mathit{f}(x)=\frac{1}{\sqrt{x}}\), then Stieltjes-integration by parts gives
\[
\sum_{p\leq x} \frac{1}{\sqrt{p}}=\frac{\pi(x)}{\sqrt{x}} -\frac{1}{\sqrt{2}}+ \frac{1}{2}\bigintss_2^x \frac{\pi(t)}{t^{3/2}} \, \mathrm{d}t,
\]
\begin{equation}\label{lim}
\implies \lim_{x\rightarrow\infty} \left( \frac{\log x}{\sqrt{x}} \sum_{p\leq x} \frac{1}{\sqrt{p}}\right)=1 +  \lim_{x\rightarrow\infty} \left( \frac{\log x}{2\sqrt{x}}  \bigintss_2^x \frac{\pi(t)}{t^{3/2}} \, \mathrm{d}t\right).
\end{equation}
Since \(\pi(x)\) is bounded by the inequality 
\[
\frac{x}{\log x + 2} < \pi(x) < \frac{x}{\log x -4},
\]
for \(x\geq 55\) \cite{rosser}, we have
$$
\lim_{x\rightarrow\infty} \left( \frac{1}{\frac{2\sqrt{x}}{\log x}}  \bigintss_2^x \bigg(\frac{1}{\sqrt{t}(\log t+2)} \bigg) \mathrm{d}t\right)< \lim_{x\rightarrow\infty} \left( \frac{\log x}{2\sqrt{x}}  \bigintss_2^x \frac{\pi(t)}{t^{3/2}} \, \mathrm{d}t\right)
$$
$$
< \lim_{x\rightarrow\infty} \left( \frac{1}{\frac{2\sqrt{x}}{\log x}} \bigintss_2^x \bigg(\frac{1}{\sqrt{t}(\log t-4)} \bigg) \mathrm{d}t\right).
$$
The bounding limits above have the indeterminate form $\frac{\infty}{\infty}$, since $\frac{1}{\sqrt{t}(\log t \,-\,4)}>\frac{1}{\sqrt{t}(\log t \,+\,2)}>\frac{1}{t}$ for $t\geq 55$  and $\int_2^x 1/t \,\mathrm{d}t$ diverges as $x\rightarrow \infty$. Therefore we can apply L'Hopital's Rule to the bounding limits, giving us
$$
1=\lim_{x\rightarrow\infty} \left(\frac{\log^2 x}{\log^2 x -4}\right) =\lim_{x\rightarrow\infty} \left(\frac{1}{\sqrt{x}(\log x+2)} \cdot \frac{\sqrt{x}\log^2 x}{\log^ x -2} \right)< \lim_{x\rightarrow\infty} \left( \frac{\log x}{2\sqrt{x}}  \bigintss_2^x \frac{\pi(t)}{t^{3/2}} \, \mathrm{d}t\right)
$$
$$
<\lim_{x\rightarrow\infty} \left(\frac{1}{\sqrt{x}(\log x-4)}\cdot \frac{\sqrt{x}\log^2 x}{\log^ x -2}\right)=\lim_{x\rightarrow\infty}\left(\frac{\log^2 x}{\log^2 x -6\log x +8}\right)=1.
$$
\p
Therefore $\lim_{x \rightarrow \infty} G(x)=2\cdot \sqrt{\pi/2}=\sqrt{2\pi}$, so our guess \eqref{guess0} becomes
 \begin{equation}\label{guess}
 \lim_{x\rightarrow\infty} \left(|\mathcal{Q}_{\mathit{f},\alpha}(x)| \frac{\log x}{\sqrt{x}}\right)=\sqrt{2\pi}.
 \end{equation}
 
 As we test our hypothesis, it should be kept in mind that $\lim_{x \rightarrow \infty}G(x)$ converges very slowly. Since the \(x\) values for which \(|\mathcal{Q}_{\mathit{f},\alpha}(x)|\frac{\log x}{\sqrt{x}}\) can actually be computed (in a reasonable amount of time) are relatively small, the largest being \(2^{27}\), we compare our computations to \(G(x)\), rather than the limit \(\sqrt{2\pi}\).
 
\section{some special cases}
In this paper we consider polynomials of the form \(\mathit{f}(z)=z^2+c\), with \(z,c\in\mathbb{Z}\), and initial argument values \(\alpha\in\mathbb{Z}\).  But for certain \(\mathit{f},\alpha\) pairs we find that we end up with a finite (over $\mathbb{Z}$) orbit, a condition which is clearly incompatible with our hypotheses outlined in sections 2 and 3, since \(m_p\) will have a fixed bound for all primes \(p\).  In this section we classify these exceptional pairs $\mathit{f},\alpha$.

\begin{prop}
Let \(\mathcal{O}_\mathit{f}(\alpha)=\{\mathit{f}^n(\alpha) : n=0,1,2,3,\dots\}\) be the orbit of \(\alpha\) under \(\mathit{f}\), where \(\mathit{f}(z)=z^2 + c, \,\,c\in\mathbb{Z}\), and \(\alpha\in\mathbb{Z}\).  Then \(\mathcal{O}_\mathit{f}(\alpha)\) is finite if and only if one of the following hold:
\begin{fleqn}
\[
\,\,\,\,\quad i) \quad\alpha = \scriptstyle\pm\frac{1}{2}\displaystyle(1\pm\scriptstyle\sqrt{1-4c}\displaystyle)
\]
\[
\,\,\quad ii)\quad\alpha = \scriptstyle\pm \frac{1}{2}\displaystyle (1 \pm \scriptstyle\sqrt{-3-4c}\displaystyle)
\]
\[
 \quad iii) \quad\alpha \in\{ 0,1,-1\} \text{ and } c\in\{0,-1,-2\}.
 \]
 \end{fleqn}
\begin{proof}
First we prove the converse, which is easier.  Assumption \((i)\) gives us the solutions to \(\alpha^2\pm\alpha + c=0\), and this equation implies that \(\alpha^2 + c = \pm\alpha\), which implies that the orbit is finite. Assumption \((ii)\) gives the solutions to \(\alpha^2 \pm \alpha + c + 1=0\), and this equation implies that \(\alpha^2 + c = \pm\alpha -1\).  With one more iteration we get 
\[
(\alpha^2 + c)^2 +c = (\pm\alpha - 1)^2 + c = \alpha^2 \mp2\alpha +1 - \alpha^2 \pm \alpha -1=\mp 2\alpha \pm \alpha = \pm \alpha,
\]
which again implies that the orbit is finite.  As for \((iii)\), testing all possible \(\alpha,\,c\) combinations will quickly convince the reader that the orbits are finite in all cases.

Now suppose that \(\mathcal{O}_\mathit{f}(\alpha)\) is finite.  First we make some simplifications.  Since the orbits of \(\alpha\) and \(-\alpha\) will be identical except for the sign of the first element \(\mathit{f}^0=\alpha\), we may consider only non-negative values of \(\alpha\).  Also, since it is obvious that \(c\in\{0,-1\}\) will have infinite orbit for \(\alpha\geq2\), and that \(c\geq 1\) will have infinite orbit for all \(\alpha\), we consider only \(c\leq-2\).  We claim that \(\mathcal{O}_\mathit{f}(\alpha)\) finite implies \(\sqrt{-c}-1<\alpha<\sqrt{-c}+1\).  If this were not true, then we would have either $\alpha=\lceil\sqrt{-c}\rceil+b$ or $\alpha=\lfloor\sqrt{-c}\rfloor-b$ for some $b\in\mathbb{N}$, giving us
\begin{align*}
\alpha=\lceil\sqrt{-c}\rceil+b \quad &\implies \quad \alpha^2 +c = (\lceil\sqrt{-c}\rceil)^2 +2b\lceil\sqrt{-c}\rceil + b^2 +c > \lceil\sqrt{-c}\rceil+b\\
\alpha=\lfloor\sqrt{-c}\rfloor-b \quad &\implies \quad \alpha^2 +c =(\lfloor\sqrt{-c}\rfloor)^2 -2b\lfloor\sqrt{-c}\rfloor + b^2 +c<-2b\lfloor\sqrt{-c}\rfloor + c.
\end{align*}
The first of these immediately implies that the iterates of \(\mathit{f}\) are unbounded since they are strictly increasing.  In the second case iterating once more gives us
\[
(\alpha^2 +c)^2 +c > 4b^2\lfloor\sqrt{-c}\rfloor^2-4bc\lfloor\sqrt{-c}\rfloor+c^2>\lfloor\sqrt{-c}\rfloor + b
\]
where the inequality reverses since \( \alpha^2 +c < -2b\lfloor\sqrt{-c}\rfloor + c<0\), and the second inequality follows since \(c\leq -2\).  Again we can conclude that the iterates of \(\mathit{f}\) are unbounded, and so we have shown that \(\mathcal{O}_\mathit{f}(\alpha)\) finite implies \(\sqrt{-c}-1<\alpha<\sqrt{-c}+1\).  For any \(c\), there are at most two integers that satisfy the preceding inequality, \(\lfloor\sqrt{-c}\rfloor\) and \(\lceil\sqrt{-c}\rceil\), so any member of \(\mathcal{O}_\mathit{f}(\alpha)\) must be one of \(\pm\lfloor\sqrt{-c}\rfloor, \pm\lceil\sqrt{-c}\rceil\), since otherwise the iterates of \(\mathit{f}\) will be unbounded.  Since we know \(\alpha \in \mathcal{O}_\mathit{f}(\alpha)\), the condition above implies that    \(\mathcal{O}_\mathit{f}(\alpha) \subset \{\alpha, -\alpha,\alpha - 1, -\alpha - 1 \}\) or  \(\mathcal{O}_\mathit{f}(\alpha) \subset \{\alpha, -\alpha,\alpha + 1, -\alpha + 1 \}\) .  However, we can rule out the latter case since 

\begin{align*}
\alpha^2 +c&= \pm\alpha +1 \\
\quad\implies\qquad\qquad (\alpha^2+c)^2 +c &= \pm3\alpha + 2 \\
\quad\implies\quad\, ((\alpha^2+c)^2+c)^2+c &= 7\alpha^2 \pm 13 \alpha +5 > 2\alpha +5 > \pm\alpha +1 > \pm \alpha,
\end{align*}
where the first inequality follows since in this case \(c\leq-2\implies \alpha\geq2\).  Therefore the iterates are unbounded in this case, and we are left with the following:

\begin{align*}
\alpha^2 +c = \pm \alpha \qquad &\text{or} \qquad \alpha^2 +c = \pm \alpha -1\\
\quad \implies \qquad \alpha^2 \pm\alpha +c =0 \qquad &\text{or} \qquad \alpha^2 \pm\alpha +c +1=0 \\
\quad \implies \quad \alpha = \scriptstyle\pm\frac{1}{2}\displaystyle(1\pm\scriptstyle\sqrt{1-4c}\displaystyle) \qquad &\text{or} \qquad \alpha = \scriptstyle\pm \frac{1}{2}\displaystyle (1 \pm \scriptstyle\sqrt{-3-4c}\displaystyle).
\end{align*}
\end{proof}
\end{prop}

This proposition is the basis for the second condition necessary for Conjectures 1 and 2; we now turn to the first condition, that \(c\notin\{0,-2\}\).  For \(c=0\) it is immediately clear that \(|\mathcal{Q}_{\mathit{f},\alpha}(x)|\) will not grow as expected, since we'll have \(p\in \mathcal{Q}_{\mathit{f},\alpha}(x) \) if and only if \(p\) divides \(\alpha\).  On the other hand, the length \(m_p\) of the orbit modulo \(p\) will grow much faster than we expect, for both \(c=0\) and \(c=-2\).  Vasiga and Shallit\cite{vas-shal} studied these cases, showing that, for a given prime \(p\), if \((p-1)/2\) is prime and \(2\) is a primitive root modulo \((p-1)/2\), then \(\sum_{0\leq\alpha<p} m_p\) is at least on the order of \(p^2\).  Heuristics by Hardy and Littlewood \cite{hardy-little}, along with Artin's conjecture, suggest that the number of primes less than \(x\) that satisfy this property is on the order of \(x/(\log x)^2\), and thus the density of these primes is on the order of \(1/\log x\).  If we sum \(p^2\), for \(p\leq x\), and multiply by \(1/\log x\), we get something on the order of \(x^3/(\log x)^2\), and dividing this by the sum \(\sum_{p\leq x} \sum_{0\leq \alpha < p} 1 \sim x^2/\log x\) gives us an average orbit length of \(\sim x/\log x\).  Note that this estimate only takes into account primes with the aforementioned property, and assumes that all other primes have orbit length \(0\), so we should expect this to be a low estimate.  Indeed, the limited experimentation we did on this question suggests that the average orbit length is closer to \(x/(\log x)^{3/4}\).

Finally, Conjecture 2 requires an additional condition, that \(\alpha^2 \neq -c\).  If we disregard this condition then we will have cases where \(0 \in \mathcal{Q}_{\mathit{f},\alpha}(x)\) for all \(p\), which is clearly incongruent with our claim.  To see that the \(\mathit{f}^0=\alpha\) is the only iterate whose square can be equal to \(-c\), suppose that the contrary were true, i.e. that we have \((\mathit{f}^l)^2 = -c\) for some \(l \in \mathbb{Z}\), then, letting \(\mathit{f}^k=\mathit{f}^{l-1}\), we have

\begin{align*}
\big((\mathit{f}^k)^2 +c\big)^2 +c &=0\\
(\mathit{f}^k)^4 +2c(\mathit{f}^k)^2  + c^2 +c &=0\\
c^2 + (2(\mathit{f}^k)^2 +1)c + (\mathit{f}^k)^4 &=0
\end{align*}
Therefore the quadratic formula gives us
\[
c=\frac{-2(\mathit{f}^k)^2-1 \pm \sqrt{(2\mathit{f}^k)^2 +1}}{2} 
\]
which is not an integer unless \(\mathit{f}^k = 0\), in which case \((\mathit{f}^l)^2=c^2= -c \implies c=-1\).  It is easy to see that this implies \(\alpha \in \{0,1\}\), and this case has already been excluded by Proposition 1(iii).
\section{results}

First we consider the first, second, third and fourth moments of $\frac{m_p}{\sqrt{p}}$, as discussed in section 2, for \(\mathit{f}(z)=z^2+c\), where \(c=\pm1,+2,\pm3\), and initial arguments \(\alpha=1,2,\dots,9\).  Of these we can exclude \(\alpha=1,2\) when \(\mathit{f}(z)=z^2-3\), and \(\alpha=1\) when \(\mathit{f}(z)=z^2-1\), because these \((\mathit{f},\alpha)\) combinations have finite orbits, as discussed above.  For the other \(42\) combinations, we find that our experimental results support our hypotheses very well.  For the first moment we expected the limit to be \(\sqrt{\pi/2}=1.25331413...\), and for all \((\mathit{f},\alpha)\) tested, \(M_1\) was between \(1.25138\) and \(1.25351\) for \(x=2^{25}\), with an average value of $1.25279$.  Table \ref{moms} gives these figures along with the standard deviation of the set of results for each moment.  It also shows the mean, standard deviation, minimum, and maximum of the set \(\{|\scriptstyle \sqrt{\pi/2} - M_1 \displaystyle |:x=2^{25}, \text{for} \,(\mathit{f},\alpha) \,\text{tested}\}\), and similarly for the second, third and fourth moments. Our complete results are depicted graphically in Figures 1-4, for the first, second, third, and fourth moments, respectively.  In each of these graphs the heavier red curve is the respective moment of \(X_n/\sqrt{n}\), for \(n=x\).  Notice that the \(y\)-axes of these graphs are not scaled equally with respect to each other (they are stretched by a factor of two for each subsequent moment graph), so if we're interested in comparing how quickly two of the moments converge, Table \ref{moms} will be more helpful.

\ctable[
caption = {Moments of $\frac{m_p}{\sqrt{p}}$ for \(x=2^{25}\) and distance from predicted limit.},
label = moms,
pos = ht
]{ccR[.]{1}{9}R[.]{1}{9}R[.]{1}{9}R[.]{1}{9}c}
{\tnote{\(\scriptstyle \sqrt{\pi/2}\displaystyle\sim 1.25331413731550\)}}
{ \FL
& & \multicolumn{1}{c}{mean}&\multicolumn{1}{c}{stand dev}& \multicolumn{1}{c}{min}&\multicolumn{1}{c}{max}&\ML
 & \(\scriptstyle M_1\)&1.25279578861576 & 0.000518157793366897 & 1.25138758183917 & 1.25350536955309 &\NN
 & \(|\scriptstyle \sqrt{\pi/2} - M_1 \displaystyle |\)&0.000544827180296359 & 0.000490240978457514 & 0.0000000515534699374598 & 0.00192655547632992 &\NN
\addlinespace
&  \(\scriptstyle M_2\)&1.99832502668792 &0.00169077550589181 & 1.9938601938444&2.00053950716997& \NN
& \(|\scriptstyle 2 - M_2\displaystyle|\)&0.00181040307318024 & 0.00154489421001528& 0.0000340794091799879&0.00613980615559995& \NN
\addlinespace
 & \(\scriptstyle M_3\)&3.75504460501883 & 0.004998322865355 & 3.74234199653338 & 3.76228591246723 & \NN
& \(|\scriptstyle 3\sqrt{\pi/2} - M_3\displaystyle |\)&0.00541926883732653 & 0.00442755795398334 & 0.0001218376421099346 & 0.0176004154131202 &\NN
\addlinespace
& \(\scriptstyle M_4\)&7.98545640132111 & 0.0149151085351382 & 7.94853101843544 & 8.00881781114142 &\NN
 & \(|\scriptstyle 8 - M_4 \displaystyle |\)&0.0162784298685631 & 0.0129995940768093 & 0.00014964867336964 & 0.0514689815645601 &\LL
}

\begin{figure}[!ht]
\labellist
\small
\pinlabel {\(\log_{\scriptscriptstyle 2} x\)} at 396 20
\pinlabel \rotatebox{90}{\(M_1\)} at 24 240
\endlabellist
\caption{The first moment of \(\frac{X}{\sqrt{n}}\) (thicker red line) and $\frac{m_p}{\sqrt{p}}$ (thin lines) for all \((\mathit{f},\alpha)\) tested.}
\label{first}
\includegraphics[scale=.55]{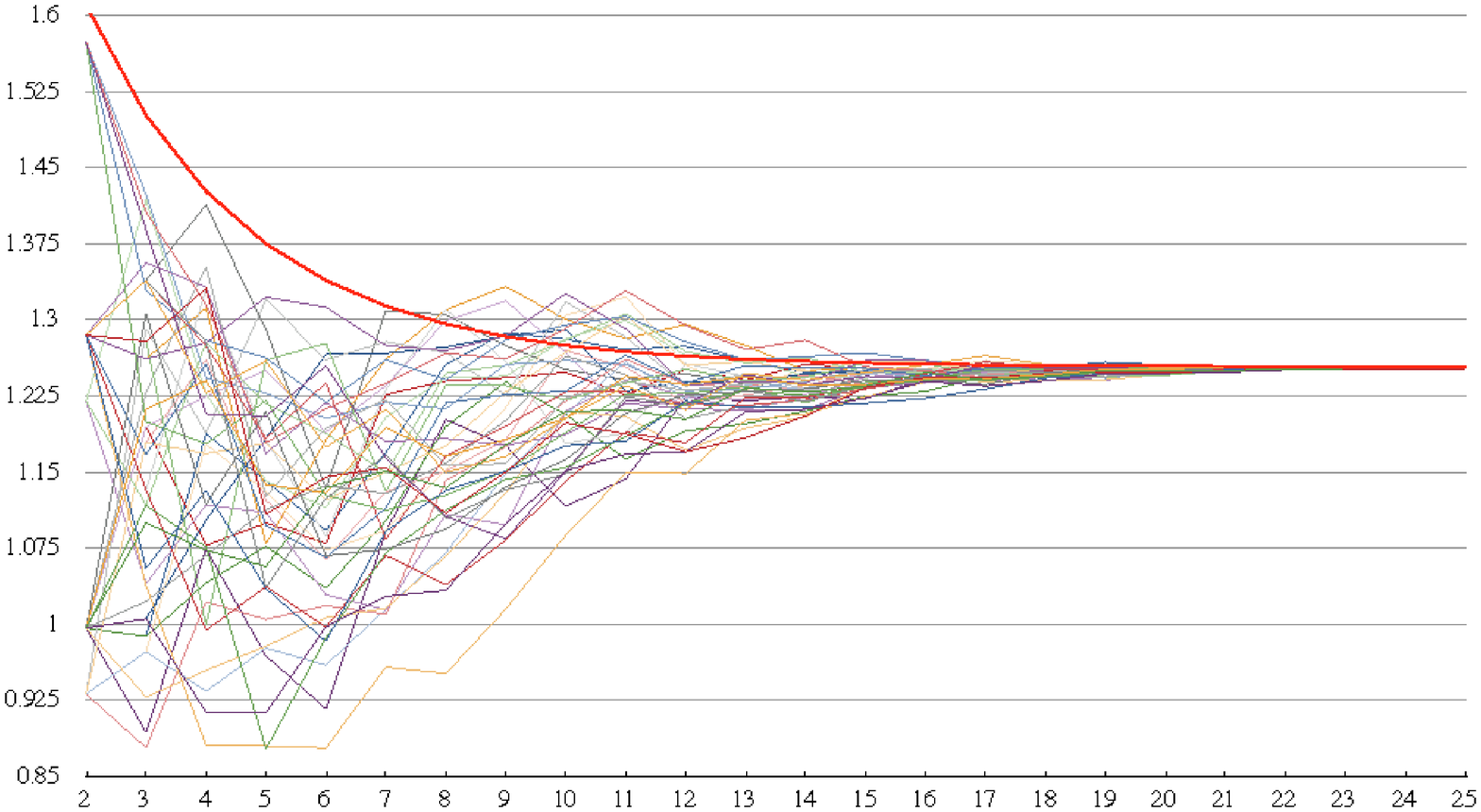}
\end{figure}

\begin{figure}[!ht]
\labellist
\small
\pinlabel {\(\log_{\scriptscriptstyle 2} x\)} at 395 18
\pinlabel \rotatebox{90}{\(M_2\)} at 19 235
\endlabellist
\caption{The second moment of \(\frac{X}{\sqrt{n}}\) (thicker red line) and $\frac{m_p}{\sqrt{p}}$ (thin lines) for all \((\mathit{f},\alpha)\) tested.}
\label{second}
\includegraphics[scale=.55]{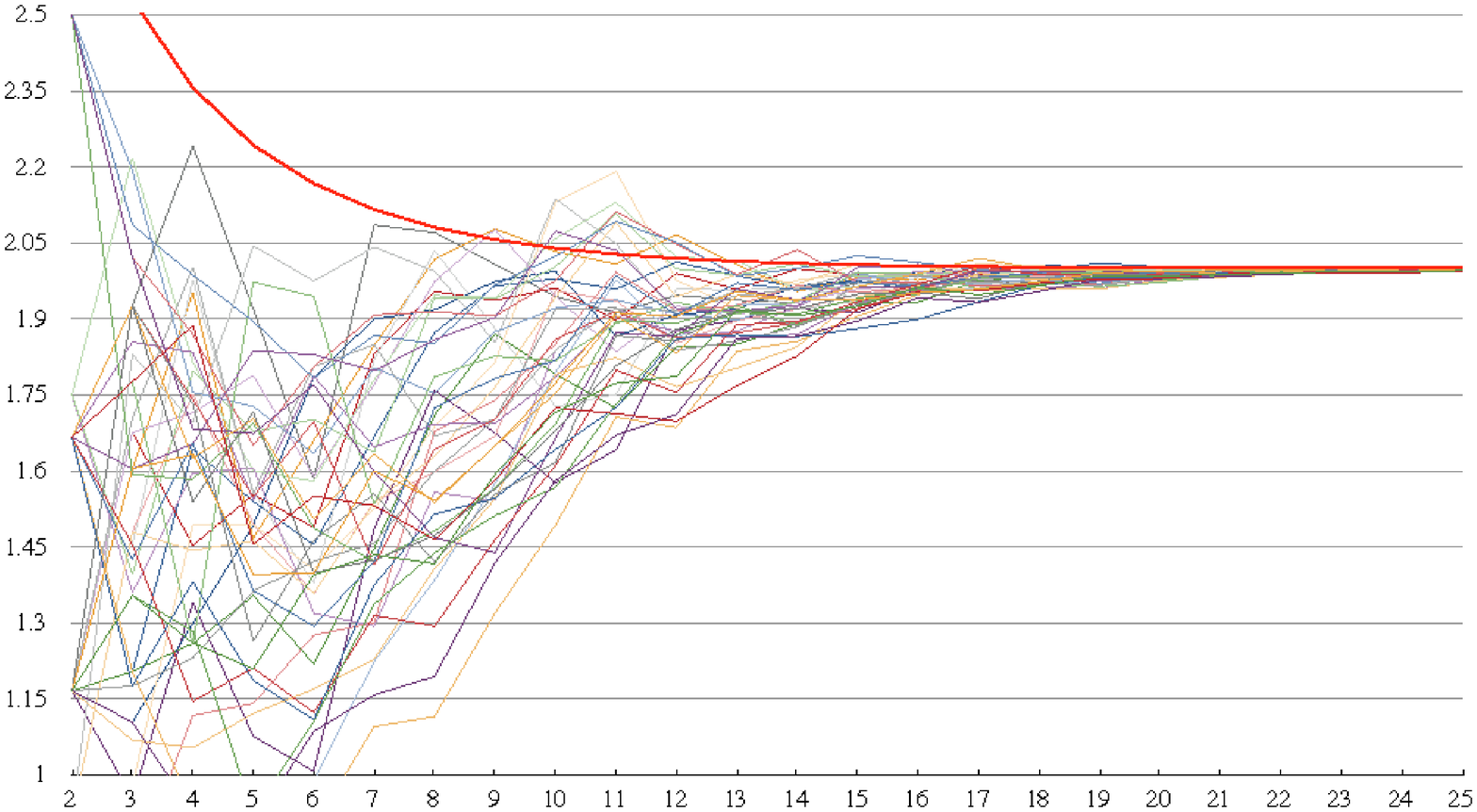}
\end{figure}

\begin{figure}[H]
\labellist
\small
\pinlabel {\(\log_{\scriptscriptstyle 2} x\)} at 393 21
\pinlabel \rotatebox{90}{\(M_3\)} at 22 224
\endlabellist
\caption{The third moment of \(\frac{X}{\sqrt{n}}\) (thicker red line) and $\frac{m_p}{\sqrt{p}}$ (thin lines) for all \((\mathit{f},\alpha)\) tested.}
\label{third}
\includegraphics[scale=.55]{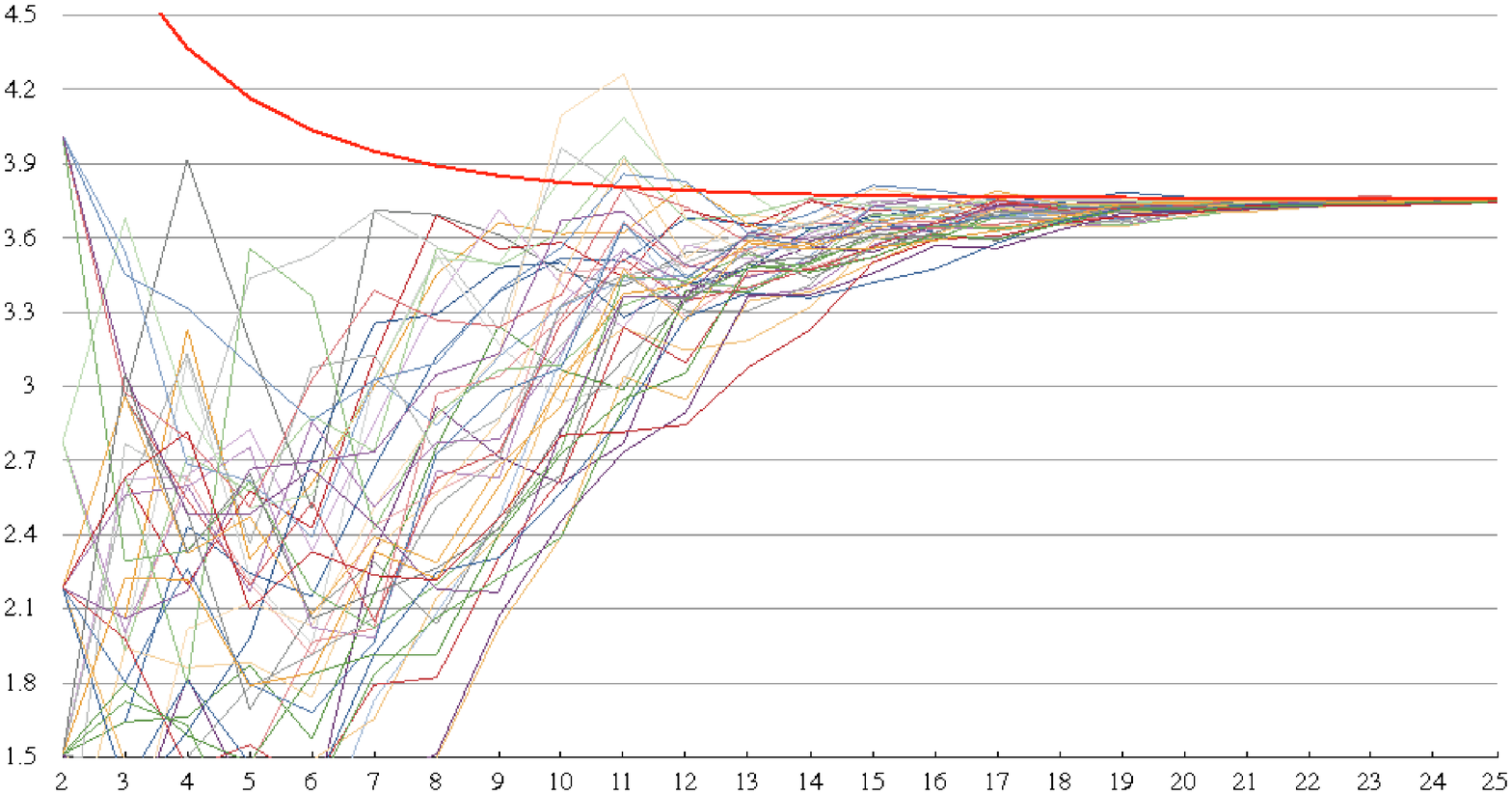}
\end{figure}

\begin{figure}[H]
\labellist
\small
\pinlabel \(\log_{\scriptscriptstyle 2}x\) at 397 26
\pinlabel \rotatebox{90}{\(M_4\)} at 23 238
\endlabellist
\caption{The fourth moment of \(\frac{X}{\sqrt{n}}\) (thicker red line) and $\frac{m_p}{\sqrt{p}}$ (thin lines) for all \((\mathit{f},\alpha)\) tested. }
\label{fourth}
\includegraphics[scale=.55]{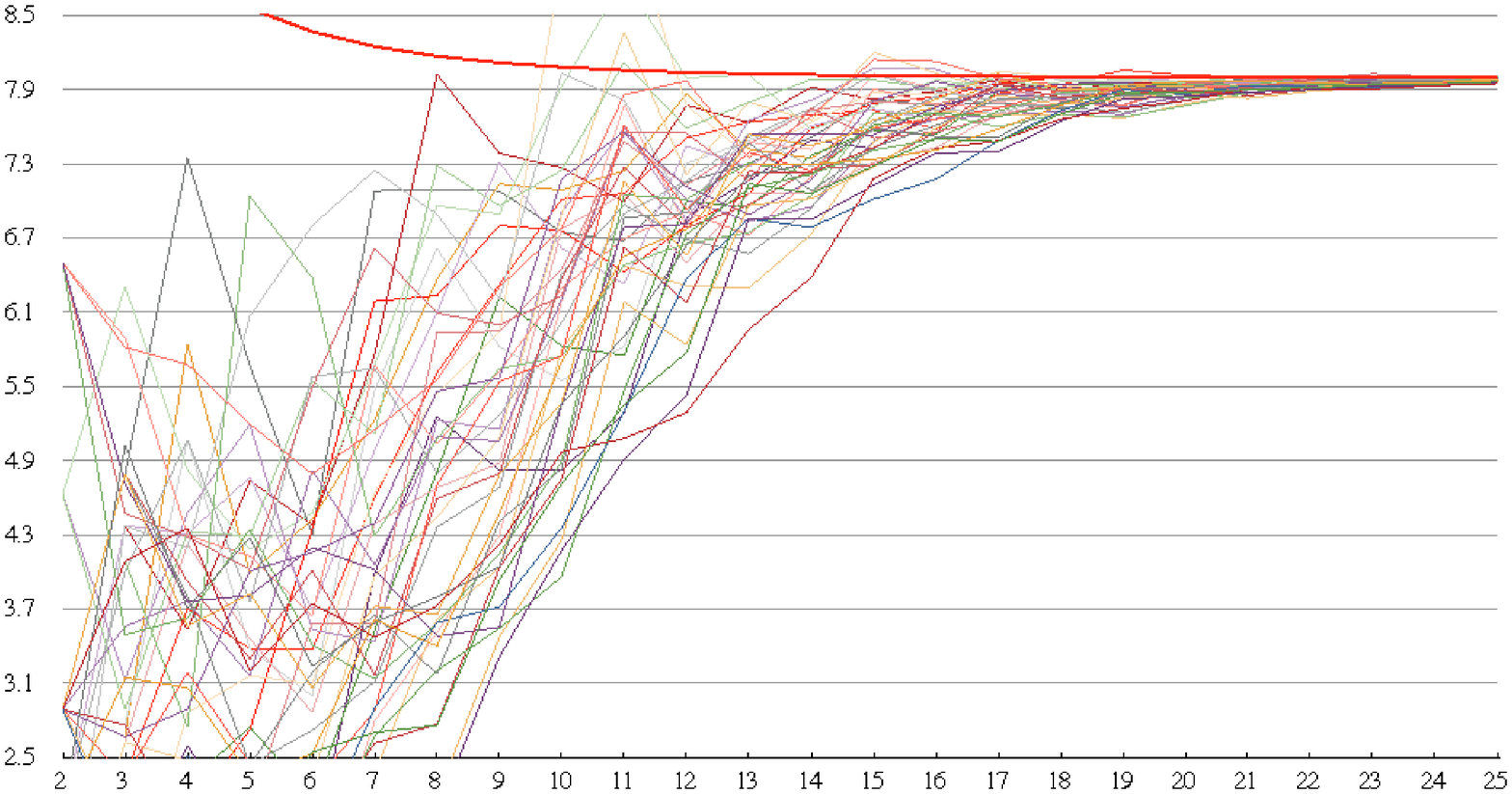}
\end{figure}

\begin{figure}[H]
\labellist
\small
\pinlabel \(t\) at 200 4
\pinlabel \rotatebox{90}{\(H(t)\) ; \(\mathit{f}(t)\)} at 0 150
\endlabellist
\caption{Histogram, H(t), of the distribution of  $\frac{m_p}{\sqrt{p}}$ (blue) for \(x=10^{8}\), \(\mathit{f}(z)=z^2+1\), \(\alpha=3\), superimposed on the graph of \(\mathit{f}(t)=te^{-t^2/2}\) (red).  \(H(t:wk \leq t < w(k+1))=\frac { | \{p \leq 10^8 : wk \leq (m_p/\sqrt{p}) < w(k+1)\}|}{w\cdot |\{p\leq 10^8\}|}\), \(k\in\mathbb{N}\). Each bar of the histogram has width \(w\approx 5.6/800\).\label{hist}}

\includegraphics[width=.89\textwidth]{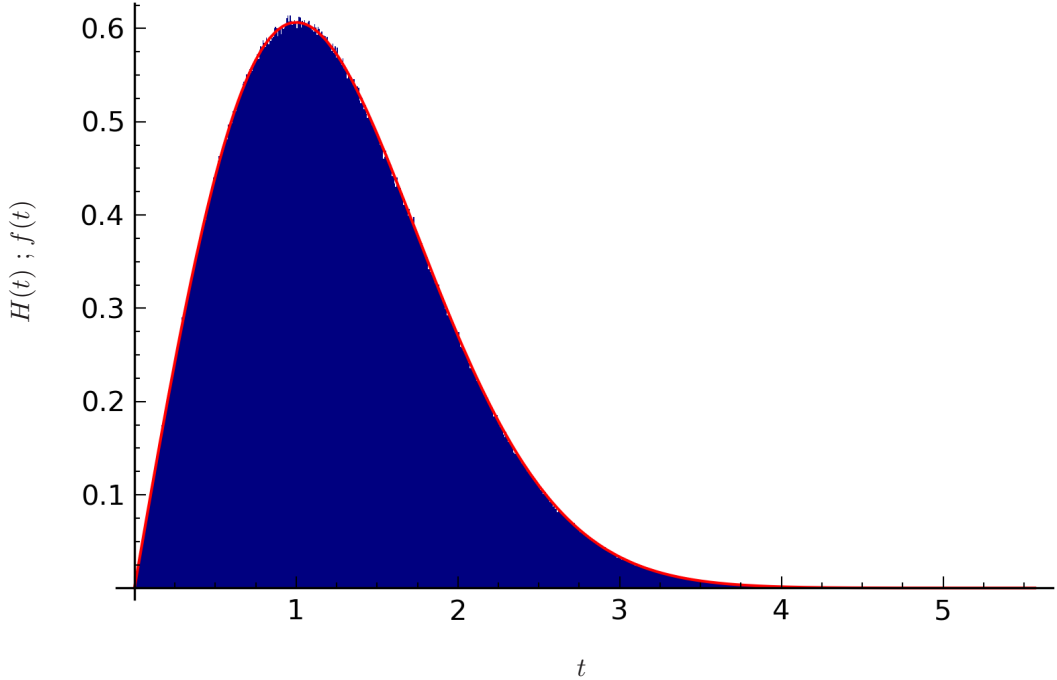}
\end{figure}

The apparent common limit of the moments of $\frac{m_p}{\sqrt{p}}$ and \(\frac{X}{\sqrt{n}}\) suggests that the limiting distributions of $\frac{m_p}{\sqrt{p}}$, as \(x \rightarrow \infty\), and the random variable \(X_n/\sqrt{n}\), as \(n \rightarrow \infty\), are the same.  For the variable \(\frac{X}{\sqrt{n}}\) we showed in section 2 that as \(n\rightarrow \infty\), the distribution \(P(X_n/\sqrt{n}<t)\) converges to the function \(F(t)=1-e^{-t^2/2}\).  As the histogram in Figure \ref{hist} shows, the density function \(\mathit{f}(t)=F'(t)\) approximates quite well the distribution of $\frac{m_p}{\sqrt{p}}$ for \(x=10^{8}\), \(\mathit{f}(z)=z^2+1\), \(\alpha=3\).  These results give considerable support to our first conjecture, stated in section 1. 

\begin{figure}[!ht]
\labellist
\small
\pinlabel \(\log_{\scriptscriptstyle 2}x\) at 352 25
\pinlabel \rotatebox{90}{\(|\mathcal{Q}_{\mathit{f},\alpha}(x)|\frac{\log x}{\sqrt{x}}\)} at 29 202
\endlabellist
\caption{Graphs of \(\mathcal{Q}_{\mathit{f},\alpha}(x)\frac{\log x}{\sqrt{x}}\) for all \(42\) \((\mathit{f},\alpha)\) combinations tested (thinner lines), the mean of these graphs (thick blue line), and our guess \(G(x)\) (thick red line).\label{quad_grph}}

\includegraphics[scale=.65]{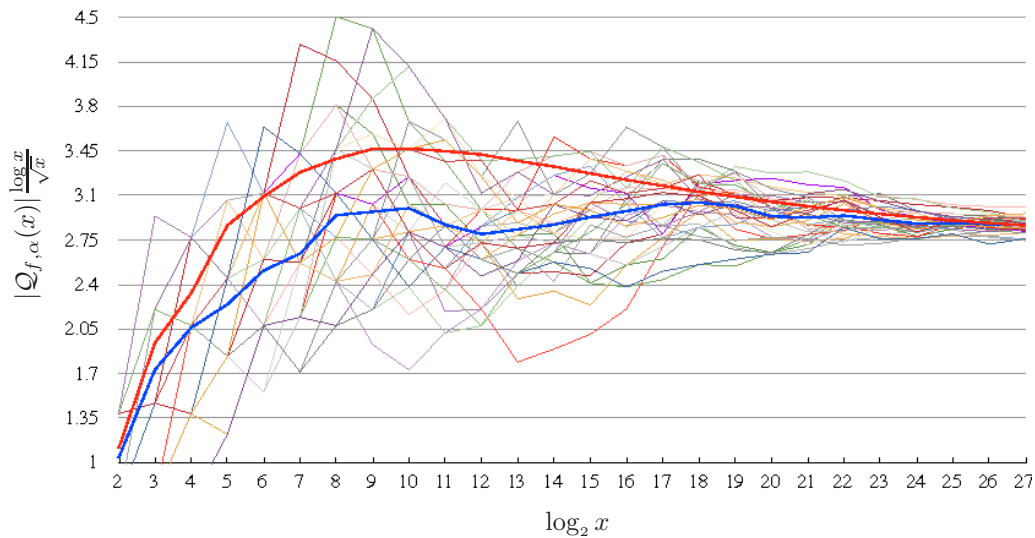}

\end{figure}

	To test the hypothesis discussed in section 3, we compute \(|\mathcal{Q}_{\mathit{f},\alpha}(x)| \frac{\log x}{\sqrt{x}}\), for \((\mathit{f},\alpha)\) as described above and \(x\in\{2,2^2,\dots,2^{27}\}\). Table \ref{quad} shows that, although our results are still fairly widely dispersed at \(x=2^{27}\), the average of the results for this \(x\) value is very close to \(G(x)\), and the standard deviation is decreasing in general as \(x\) increases, as is the error of the mean from \(G(x)\). As we mentioned earlier, \(\lim_{x\rightarrow\infty} G(x)\) converges very slowly, and as the table shows, even for \(x\) as large as \(2^{27}\) we still have \(|G(x)-\sqrt{2\pi}|\sim 0.36\), so we are not too surprised to see such a wide range in our results for this \(x\) value.  That is, intuitively, it seems we should not expect our results to be very tightly grouped until we are close to the limiting value, \(\sqrt{2\pi}\).  Figure \ref{quad_grph} gives a graphical representation of all \((\mathit{f},\alpha)\) tested, for \(x\) from \(4\) to \(2^{27}\). On this graph the red and blue lines are $G(x)$ and the mean from Table \ref{quad}, respectively.  From this data, it seems reasonable to suppose that \(|\mathcal{Q}_{\mathit{f},\alpha}(x)| \frac{\log x}{\sqrt{x}}\) will eventually converge to \(\sqrt{2\pi}\), independent of \(\mathit{f}\), \(\alpha\), and so we make our second conjecture as stated in section 1.

In his paper \emph{Variation of Periods Modulo p in Arithmetic Dynamics}, Silverman carries out computations that lead to a conjecture (in a more general setting) that under certain restrictions the set $\{p: m_p \leq p^{1/2-\epsilon}\}$ will have density $0$ for $\epsilon>0$ \cite{silv}. This conjecture agrees with our own results, and in fact if Conjecture 1 were proven, a less general version of Silverman's conjecture would readily follow. Computations of a similar nature to ours were also carried out in \emph{Periods of Rational Maps Modulo Primes} by Benedetto, et al, with results that are compatible with our own \cite{hutz}.

\ctable
[
caption = {A comparison of our experimental results to our guess, \(G(x)\).  },
label = quad,
pos = ht
]
{ccR[.]{1}{5}R[.]{1}{5}R[.]{1}{5}R[.]{1}{5}R[.]{1}{5}R[.]{1}{5}c}
{\tnote{\(\frac{\log x}{\sqrt{x}} \sum_{p\leq x} \frac{\sqrt{\scriptstyle \pi/2}}{\sqrt{p}}\)}
}
{ 
\FL
&\multicolumn{1}{c}{} &\multicolumn{1}{c}{}& \multicolumn{4}{c}{\(\scriptstyle|\mathcal{Q}_{\mathit{f},\alpha}(x)| \frac{\log x}{\sqrt{x}}\) for all \((\mathit{f},\alpha)\) tested} &\multicolumn{1}{c}{} &\NN
\cmidrule(rl){4-7}
& \(x\) & \multicolumn{1}{c}{\(\quad\)\(G(x)\)\tmark\(\quad\)} & \multicolumn{1}{c}{\(\,\,\)mean\(\,\,\)} & \multicolumn{1}{c}{\(\,\,\)stand dev\(\,\,\)} & \multicolumn{1}{c}{\(\,\,\)minimum\(\,\,\)}& \multicolumn{1}{c}{\(\,\,\)maximum\(\,\,\)}& \multicolumn{1}{c}{\(\, | G(x) - \text{mean}  |\)} & \ML
&\(2^{10}\)&3.46924936385524&3.00157484438905&0.516867879926488&1.73286795139986&4.11556138457468&0.467674519466191&\NN
&\(2^{11}\)&3.45002662400205&2.87221470906186&0.469325204349483&2.02178242090388&3.70660110499044&0.577811914940192&\NN
&\(2^{12}\)&3.42304422621934&2.79734397868835&0.34646305811633&2.07944154167984&3.37909250522973&0.625700247530992&\NN
&\(2^{13}\)&3.37313224613506&2.83501741343917&0.375185246911875&1.79203441852844&3.6836263047529&0.538114832695891&\NN
&\(2^{14}\)&3.32853583397385&2.87186761528873&0.309152297102417&1.8953243218436&3.56320972506597&0.456668218685119&\NN
&\(2^{15}\)&3.27415308950979&2.93202212555325&0.350711316671737&2.01029502078511&3.44622003563161&0.342130963956545&\NN
&\(2^{16}\)&3.22424881093474&2.97785105695917&0.333974632293458&2.20940663803483&3.63902269793971&0.24639775397557&\NN
&\(2^{17}\)&3.17737202507932&3.03080464959046&0.288609036529303&2.44107252523906&3.48259680267439&0.146567375488857&\NN
&\(2^{18}\)&3.13140372042091&3.04547674449706&0.188494007059308&2.55868783448886&3.38721532375192&0.0859269759238508&\NN
&\(2^{19}\)&3.09044792152051&3.01797249429906&0.210036095568691&2.54637369299447&3.3284741844142&0.0724754272214532&\NN
&\(2^{20}\)&3.05137188874281&2.93775269885758&0.176019194496016&2.63991601971073&3.27620347061537&0.113619189885228&\NN
&\(2^{21}\)&3.01654146201451&2.92737246360041&0.157861799402495&2.65358942743634&3.28683235898365&0.0891689984141029&\NN
&\(2^{22}\)&2.98474862867106&2.94397374637656&0.129884507838308&2.71031378023635&3.216636134786&0.0407748822945022&\NN
&\(2^{23}\)&2.95589166360345&2.91037498614292&0.0940209747359812&2.72466771567124&3.1154786405453&0.0455166774605287&\NN
&\(2^{24}\)&2.92985593883531&2.88060286854077&0.0842753933529274&2.76175829754353&3.07448681064773&0.0492530702945428&\NN
&\(2^{25}\)&2.90646078014942&2.88289012777343&0.0544511570308805&2.77612169536991&3.02740857296805&0.023570652375994&\NN
&\(2^{26}\)&2.88509497148514&2.87750846148469&0.0586916386831027&2.71911350030205&3.01390412250308&0.00758651000044619&\NN
&\(2^{27}\)&2.86577988450825&2.86820893411965&0.0541783737649122&2.74620659089394&3.00790392484972&0.002429049611401&\LL}

\emph{}

\textsc{B.S. program in Mathematics, City College of New York, Convent Ave at 138th St, New York, NY 10031.}

\emph{E-mail address:} \texttt{wworden00@ccny.cuny.edu}

\end{document}